\documentclass[12pt]{amsart}

\usepackage{upref, amssymb}
\usepackage{graphicx}

\usepackage{hyperref}

\newcommand{\dmat}[1]{\setbox0\hbox{\kern1pt $#1$\kern1pt
}\mathord{\hbox to 0pt{\kern1pt $#1$\hss}\vrule width\wd0
height2.5pt depth -2pt}}

\newcommand{\dmatCAP}[1]{\setbox0\hbox{\kern1pt $#1$\kern1pt
}\mathord{\hbox to 0pt{\kern1pt $#1$\hss}\vrule width\wd0
height3.5pt depth -3pt}}

\newcommand{\umat}[1]{\overline{#1}}
\newcommand{\vmat}[1]{\underline{#1}}

\newcommand{\<}{\mathbin{\ <\ }}

\newcommand{\integer}{{\mathord{\mathbb{Z}}}}

\newcommand{\rational}{{\mathord{\mathbb{Q}}}}
\newcommand{\real}{{\mathord{\mathbb{R}}}}
\newcommand{\complex}{{\mathord{\mathbb{C}}}}

\newcommand{\tuple}{\overrightarrow}
\newcommand{\bigtimes}{\mathop{\text{\Large$\times$}}\limits}

\renewcommand{\pmod}[1]{\ ({\rm mod} \  #1)}
\newcommand{\bigpmod}[1]{\ \bigl({\rm mod} \  #1\bigr)}

\newcommand{\Qrank}{\mathop{\mbox{\upshape$\rational$-rank}}}
\newcommand{\Rrank}{\mathop{\mbox{\upshape$\real$-rank}}}

\newcommand{\Frank}{\mathop{\mbox{\upshape$\F$-rank}}}

\newcommand{\cjg}[1]{\overline{#1}}
\newcommand{\iso}{\cong}

\newcommand{\Id}{\mathord{\mathrm{Id}}}
\newcommand{\SU}{\operatorname{SU}}

\newcommand{\SUJ}{\SU_{2,1}}

\newcommand{\ints}{\mathcal{O}}
\DeclareMathOperator{\SL}{SL}
\DeclareMathOperator{\PSL}{PSL}
\newcommand{\transpose}{T}

\DeclareMathOperator{\Ad}{Ad}

\newcommand{\U}{\mathcal{U}}
\newcommand{\V}{\mathcal{V}}
\newcommand{\T}{\mathcal{T}}

\newcommand{\scrS}{\mathcal{S}}

\newcommand{\G}{\mathbf{G}}
\newcommand{\GU}{\mathbf{U}}
\newcommand{\GV}{\mathbf{V}}
\newcommand{\GT}{\mathbf{T}}

\newcommand{\F}{\mathbb{F}}
\newcommand{\K}{\mathbb{K}}

\makeatletter
\renewcommand{\p@enumii}{} 
\makeatother

\newcommand{\bigset}[2]{\left\{\, #1 
 \mathrel{\left| \vphantom {\left\{ #1 \mid #2 \right\} }
 \right.} #2 \,\right\} }

\renewcommand{\see}[1]{{\upshape(}see~\ref{#1}{\upshape)}}
\newcommand{\seeand}[2]{{\upshape(}see~\ref{#1} and~\ref{#2}{\upshape)}}
\newcommand{\cf}[1]{{\upshape(}cf.~\ref{#1}{\upshape)}}
\newcommand{\fullsee}[2]{{\upshape(}see~\ref{#1}\pref{#1-#2}{\upshape)}}
\newcommand{\fullcf}[2]{{\upshape(}cf.~\ref{#1}\pref{#1-#2}{\upshape)}}
\newcommand{\fullfullcf}[3]{{\upshape(}cf.~\ref{#1}{\upshape(}\ref{#1-#2}\ref{#1-#2-#3}{\upshape)}{\upshape)}}

\newcommand{\cfStep}[1]{{\upshape(}cf.\ Step~\ref{#1}{\upshape)}}

\newcommand{\pref}[1]{{\upshape(}\ref{#1}{\upshape)}}
\newcommand{\fullref}[2]{{\ref{#1}\pref{#1-#2}}}

 \makeatletter
 \def\swappedhead#1#2#3{%
  \@ifnotempty{#2}{{\normalfont(#2)}}{}%
  \thmname{\@ifnotempty{#2}{~}#1}%
  \thmnote{ {\the\thm@notefont{#3}}}} 
 \makeatother

\swapnumbers

\newtheorem{thm}[equation]{Theorem}
\newtheorem{prop}[equation]{Proposition}
\newtheorem{cor}[equation]{Corollary}
\newtheorem{lem}[equation]{Lemma}
\newtheorem{conj}[equation]{Conjecture}

\theoremstyle{definition}
\newtheorem{defn}[equation]{Definition}
\newtheorem{notation}[equation]{Notation}
\newtheorem{rem}[equation]{Remark}
\newtheorem{eg}[equation]{Example}
\newtheorem{ack}[equation]{Acknowledgments}
\newtheorem{assump}[equation]{Assumption}

 \newtheorem{noheading}[equation]{}
 
 \newenvironment{proofnum}[1]{\refstepcounter{equation}
\begin{proof}[{\rm(\theequation) {\bf Proof of #1}}]}{\end{proof}}

\numberwithin{equation}{section}

 \newcounter{step}

\newenvironment{step}[1][\unskip]{\refstepcounter{step}
 \setcounter{substep}{0}
 \em
 \medskip \noindent Step \thestep\ #1.\ }{\unskip\upshape}
 \renewcommand{\thestep}{\arabic{step}}

 \newcounter{substep}

 \renewcommand{\thesubstep}{\thestep.\arabic{substep}}

 \newcounter{case}

 \renewcommand{\thecase}{\arabic{case}}

 \newcounter{subcase}
 
 \renewcommand{\thesubcase}{\arabic{case}.\arabic{subcase}}

\newcommand{\thmrefer}[1]{\renewcommand\theequation
  {\protect\ref{#1}$'$}\addtocounter{equation}{-1}}

\begin{document}

\title
 [Bounded generation and actions on the line]
 {Bounded generation and lattices \\
 that cannot act on the line}
 
 \dedicatory{To Professor G.\,A.\,Margulis on his 60th birthday}

\author{Lucy Lifschitz}

\address{Department of Mathematics,
 University of Oklahoma,
 Norman, Oklahoma, 73019, USA}
 \email{llifschitz@math.ou.edu,
 http://www.math.ou.edu/$\sim$llifschitz/}

\author{Dave Witte Morris}

\address{Department of Mathematics and Computer Science,
 University of Lethbridge,
 Lethbridge, Alberta, T1K~3M4, Canada}
\email{Dave.Morris@uleth.ca,
 http://people.uleth.ca/$\sim$dave.morris/}
 
 \begin{abstract}
 Let $\Gamma$ be an irreducible lattice in a connected, semisimple Lie group with finite center. 
 Assume that $\Rrank G \ge 2$, that $G/\Gamma$ is not compact, and that $G$ has more than one noncompact simple factor.
 We show that $\Gamma$ has no orientation-preserving actions on the real line. 
 (In algebraic terms, this means that $\Gamma$ is not right orderable.)
Under the additional assumption that no simple factor of~$G$ is isogenous to $\SL(2,\real)$,
applying a theorem of \'E.~Ghys yields the conclusion that any orientation-preserving action of~$\Gamma$ on the circle must factor through a finite, abelian quotient of~$\Gamma$. 

The proof relies on the fact, proved by D.~Carter, G.~Keller, and E.~Paige, that $\SL(2,\ints)$ is boundedly generated by unipotents whenever $\ints$ is a ring of integers with infinitely many units.
The assumption that $G$ has more than one noncompact simple factor can be eliminated if all noncocompact lattices in $\SL(3,\real)$ and $\SL(3,\complex)$ are virtually boundedly generated by unipotents.
 \end{abstract}

 \date{June 16, 2007} 

\maketitle

\section{Introduction} \label{IntroSect}

It is known that if $\Gamma$ is a finite-index subgroup of $\SL(3,\integer)$, then $\Gamma$ has no nontrivial actions by orientation-preserving homeomorphisms of the real line~$\real$ \cite{Witte-ActOnCircle}. (More generally, the same is true if $\Gamma$ is any finite-index subgroup of the integer points of any connected, almost-simple, algebraic group over~$\rational$, with $\Qrank G \ge 2$.) It is conjectured that the same conclusion is true much more generally:

\begin{defn}
A subgroup~$\Gamma$ of a Lie group~$G$ is an \emph{irreducible lattice} in~$G$ if
\begin{enumerate}
\item $\Gamma$ is discrete,
\item $G/\Gamma$ has finite volume,
and
\item $\Gamma N$ is dense in~$G$, for every noncompact, closed, normal subgroup~$N$ of~$G$.
\end{enumerate}
\end{defn}

\begin{conj}[\cite{GhysAOC}] \label{Rrank2LieConj}
Suppose
\begin{itemize}
\item $G$ is a connected, semisimple Lie group with finite center,
\item $\Rrank G \ge 2$,
and
\item $\Gamma$ is any irreducible lattice in~$G$.
\end{itemize}
Then\/ $\Gamma$ has no nontrivial, orientation-preserving action on\/~$\real$.
\end{conj}

In this paper, we prove the conjecture in the special case where $G$ is a direct product of copies of $\SL(2,\real)$ and/or $\SL(2,\complex)$.

\begin{eg}
The following theorem implies that no finite-index subgroup of $\SL \bigl( 2, \integer[\sqrt{3}] \bigr)$ has a nontrivial, orientation-preserving action on~$\real$. (Such subgroups are noncocompact, irreducible lattices in $\SL(2,\real) \times \SL(2,\real)$). Furthermore, $\sqrt{3}$ can be replaced with any irrational algebraic integer~$\alpha$, such that {either} $\alpha$~is real {or} $\alpha$ is \emph{not} a root of any quadratic polynomial with rational coefficients.
\end{eg}

\begin{thm} \label{SL(2,O)NoActR}
Let
\begin{itemize}
\item $\F$ be an algebraic number field that is neither\/ $\rational$ nor an imaginary quadratic extension of\/~$\rational$,
\item $\ints$ be the ring of integers of\/~$\F$,
and
\item $\Gamma$ be a finite-index subgroup of\/ $\SL(2,\ints)$.
\end{itemize}
Then\/ $\Gamma$ has no nontrivial, orientation-preserving action on\/~$\real$. 
\end{thm}

\begin{rem} \ 
\begin{enumerate}
\item See \cite{GhysSurvey} for a very nice introduction to this subject.
\item A version of Conjecture~\ref{Rrank2LieConj} circulated informally in 1990, but it apparently first appeared in print in~\cite{GhysAOC}.
\item The conclusion of the conjecture is equivalent to the purely algebraic statement that $\Gamma$ is not right orderable \see{Act<>RO}.
That is, there does not exist a total order~$\prec$ on~$\Gamma$, such that $a \prec b$ implies $ac \prec bc$, for all $c \in \Gamma$. 
\item  Theorem~\ref{SL(2,O)NoActR} was announced in \cite{LifschitzMorris-CR}. It provides the first known examples of arithmetic groups of $\Qrank$~1 that have no right-orderable subgroups of finite index.
\item For an algebraic number field~$\F$ with ring of integers~$\ints$, the Dirichlet Units Theorem (cf.\ \cite[Prop.~4.7, p.~207]{PlatonovRapinchukBook}) implies that the following two conditions are equivalent:
\begin{enumerate}
\item $\F$ is neither\/ $\rational$ nor an imaginary quadratic extension of\/~$\rational$.
\item The group $\ints^\times$ of units of~$\ints$ is infinite.
\end{enumerate}
\end{enumerate}
\end{rem}

The above theorem considers only a very restricted class of lattices. However, because subgroups of right orderable groups are right orderable, it has more general consequences. For example, it implies that the conclusion of the conjecture holds when $G$ has more than one simple factor and $\Gamma$ is not cocompact:

\begin{cor}[{\cf{NonsimpleContainsSL2}}] \label{NotSimple->NotRO}
Assume
\begin{itemize}
\item $G$ and\/ $\Gamma$ are as in Conjecture~\ref{Rrank2LieConj},
\item the adjoint group of~$G$ is \textbf{not} simple,
and
\item $G/\Gamma$ is \textbf{not} compact.
\end{itemize}
Then\/ $\Gamma$ has no nontrivial, orientation-preserving action on\/~$\real$.
\end{cor}

The proof of Theorem~\ref{SL(2,O)NoActR} is based on the following simple lemma.

\begin{lem}[{\see{BddGen+BddOrbitsPf}}] \label{BddGen+BddOrbits}
Suppose
\begin{itemize}
\item $\Gamma$ is a group,
\item $U_1,U_2,\ldots,U_r$ are subgroups of~$\Gamma$,
\item the product $U_1 U_2 \cdots U_r$ is a finite-index subgroup of\/~$\Gamma$,
and
\item for every orientation-preserving action of\/~$\Gamma$ on\/~$\real$, and for each $i \in \{1,\ldots,r\}$, the $U_i$-orbit of each point in\/~$\real$ is a bounded set.
\end{itemize}
Then\/ $\Gamma$ has no nontrivial, orientation-preserving action on\/~$\real$.
\end{lem}

Thus, Theorem~\ref{SL(2,O)NoActR} is a consequence of the following two theorems. Before stating these results, we provide an important definition.

\begin{defn} \ 
\begin{itemize}
\item A subgroup $U$ of $\SL(\ell,\complex)$ is \emph{unipotent} if it is conjugate to a subgroup of 
$$ \begin{bmatrix}
1 &  & \vbox to 0pt{\vskip-5pt\hbox to 0pt{\hss\mbox{\huge{$*$}}}\vss}  \\
 & \ddots  \\
\vbox to 0pt{\vss\hbox to 0pt{\mbox{\huge{$0$}}\hss}} & & 1
\end{bmatrix} .$$
\item A matrix group $\Gamma$ is \emph{virtually boundedly generated by unipotents} if there are unipotent subgroups $U_1,\ldots,U_r$ of~$\Gamma$, such that the product $U_1 U_2 \cdots U_r$ is a finite-index subgroup of~$\Gamma$.
\end{itemize}
\end{defn}

\begin{thm}[{(D.~Carter, G.~Keller, and E.~Paige \cite{CarterKellerPaige, MorrisBddGenSL2})}] \label{SL2BddGen}
If\/ $\F$, $\ints$, and\/~$\Gamma$ are as as described in Theorem~\ref{SL(2,O)NoActR}, then\/ $\Gamma$ is virtually boundedly generated by unipotents.
\end{thm}

\begin{thm}[(see \S\ref{ArchimedeanUBddOrbitsPf})] \label{ArchimedeanUBddOrbits}
Suppose
\begin{itemize}
\item $\F$, $\ints$, and\/~$\Gamma$ are as described in Theorem~\ref{SL(2,O)NoActR}, 
\item no proper subfield of\/~$\F$ contains a finite-index subgroup of the group~$\ints^\times$ of units of~$\ints$,
and
\item $U$ is any unipotent subgroup of\/~$\Gamma$.
\end{itemize}
Then, for every orientation-preserving action of\/~$\Gamma$ on\/~$\real$, the $U$-orbit of each point in\/~$\real$ is a bounded set.
\end{thm}

The following theorem shows that our methods will yield more general results if one can generalize the Carter-Keller-Paige Theorem \pref{SL2BddGen} to establish the bounded generation of additional groups.

\begin{conj} \label{SU2conj}
If\/ $\Gamma$ is any noncocompact lattice in either\/ $\SL(3,\real)$ or\/ $\SL(3,\complex)$, then\/ $\Gamma$ is virtually boundedly generated by unipotents.
\end{conj}

\begin{thm}[(see \S\ref{AllLattSect})] \label{AllNotRO}
Assume
\begin{itemize}
\item Conjecture~\ref{SU2conj} is true,
\item $G$ is a connected, semisimple Lie group with finite center,
\item $\Rrank G \ge 2$,
and
\item $\Gamma$ is a noncocompact, irreducible lattice in~$G$.
\end{itemize}
Then\/ $\Gamma$ has no nontrivial orientation-preserving action on\/~$\real$.
\end{thm}

\begin{rem} \ 
\begin{enumerate}
 \item The Margulis Arithmeticity Theorem provides a concrete description of the noncocompact lattices in $\SL(3,\real)$ and $\SL(3,\complex)$ \cf{SL3existsKF}.
 \item Conjecture~\ref{SU2conj} is only a very special case of a much more general conjecture: it is believed that $\SL(3,\real)$ and $\SL(3,\complex)$ can be replaced by any simple Lie groups of real rank $\ge 2$ (cf.~\cite[p.~578]{PlatonovRapinchukBook}).
 \end{enumerate}
\end{rem}

A beautiful theorem of \'E.~Ghys \cite[Thm.~3.1]{GhysAOC} implies that if $\Gamma$ is a higher-rank lattice, then (up to finite covers) any action of~$\Gamma$ on the circle~$S^1$ must be semiconjugate to an action obtained from projecting to a $\PSL(2,\real)$ factor of~$G$.
 In particular, if there are no nontrivial homomorphisms from~$G$ to $\PSL(2,\real)$, then every action of~$\Gamma$ on~$S^1$ has a finite orbit. 
 (In most cases, this conclusion was also proved by M.~Burger and N.~Monod \cite{BurgerMonod1, BurgerMonod2}.)
Combining this with Theorem~\ref{AllNotRO} yields the following conclusion:

\begin{cor} \label{NoSL2->NoActS1}
Assume 
\begin{itemize}
\item $\Gamma$ and~$G$ are as in Theorem.~\ref{AllNotRO},
\item Conjecture~\ref{SU2conj} is true,
 and
 \item no simple factor of~$G$ is isogenous to\/ $\SL(2,\real)$.
\end{itemize}
Then any action of\/~$\Gamma$ on the circle~$S^1$ factors through a finite quotient of\/~$\Gamma$.
\end{cor}

Regrettably, our methods do not apply to cocompact lattices, because these do not have any unipotent subgroups.

Here is an outline of the paper: 
\begin{enumerate}
\item[\S\ref{IntroSect}.] Introduction
\item[\S\ref{SarithSect}.] The $S$-arithmetic case
\item[\S\ref{PrelimArithSect}.] Preliminaries on arithmetic groups
\item[\S\ref{PrelimBddGenSect}.] Preliminaries on bounded generation
\item[\S\ref{PrelimUnbddOrbitSect}.] Preliminaries on unbounded orbits of unipotent subgroups
\item[\S\ref{ArchimedeanUBddOrbitsPf}.] Proof of Theorem~\ref{ArchimedeanUBddOrbits}
\item[\S\ref{LattSL3NotRO}.] Lattices in $\SL(3,\real)$ or $\SL(3,\complex)$
\item[\S\ref{AllLattSect}.] Proof of Theorem~\ref{AllNotRO}
\end{enumerate}

\begin{ack}
We thank V.~Chernousov for very helpful conversations. 
The work of D.W.M.\ was partially supported by a grant from the
National Sciences and Engineering Research Council of Canada.
\end{ack}

\section{The $S$-arithmetic case} \label{SarithSect}

As an easy introduction to the methods that prove Theorem~\ref{ArchimedeanUBddOrbits}, let us first consider the situation where the ring~$\ints$ of integers is replaced with a ring $\integer[1/r]$ of $S$-integers (with $r \neq \pm 1$). (Thus, $\Gamma$ is an $S$-arithmetic group, rather than an arithmetic group.) B.~Liehl \cite{LiehlBddGenSL2} proved bounded generation by unipotents in this setting, so we conclude that $\Gamma$ has no nontrivial actions on~$\real$ \see{SL(2,Z[1/r])NotRO}. This yields analogues of Corollaries~\ref{NotSimple->NotRO} and \ref{NoSL2->NoActS1} in which some of the simple factors of~$G$ are $p$-adic, rather than real \seeand{SarithThm}{SarithNoActS1}. All of these results appeared in \cite{LifschitzMorris-CR}.

\begin{prop}[{\cite[Thm.~1.4(i)]{LifschitzMorris-CR}}] \label{SL(2,Z[1/r])UBddOrb}
 Let\/ $\Gamma$ be a finite-index subgroup of\/ $\SL \bigl( 2, \integer[1/r] \bigr)$, for some natural number~$r>1$.

For each action of\/~$\Gamma$ on\/~$\real$, every orbit of every unipotent
subgroup of\/~$\Gamma$ is bounded.
 \end{prop}

\begin{proof}
 Suppose $\Gamma$ acts on~$\real$, and, for some unipotent subgroup~$U_1$ of~$\Gamma$, the $U_1$-orbit of some point~$x$ is not bounded.
(This will lead to a contradiction.) 
 We begin by establishing notation.
 \begin{itemize}
  
 \item For $u,v,w \in \rational$, with $w \neq 0$, let 
 $$ \mbox{$\umat{u} = \begin{bmatrix}
 1 & u \\
 0 & 1 \\
 \end{bmatrix}$,
 \qquad
 $\vmat{v} = \begin{bmatrix}
 1 & 0 \\
 -v & 1 \\
 \end{bmatrix}$
  \quad and \quad
 $\dmat{w} = \begin{bmatrix} w & 0 \\ 0 & 1/w \end{bmatrix}$%
 .}$$
 Note that $\umat{u}$, $\vmat{v}$, and $\dmat{w}$ each belong to $\SL(2, \rational)$.
 (The minus sign in the definition of~$\vmat{v}$ ensures that $\umat{u}$ is conjugate to~$\vmat{v}$ when $u = v$ \see{upperlowerconjugate}.)

 \item Let
 $$ \mbox{$\U = \bigset{ \umat{u} }{ u \in \rational }$
 \quad and \quad
 $\V = \bigset{ \vmat{v} }{ v \in \rational }$},$$
 so $\U$ and~$\V$ are opposite maximal unipotent subgroups of
$\SL(2,\rational)$. 

\item Let
 $$ \mbox{$U = \U \cap \Gamma$ \quad and \quad $V = \V \cap \Gamma$.}$$

 \item Fix some $\omega \in \{\, r^n \mid n \in \integer^+ \,\}$, such that $\dmat{\omega} \in \Gamma$ (this is possible because $\Gamma$ has finite index in $\SL \bigl( 2, \integer[1/r] \bigr)$.
 Note that $\omega > 1$ (because $n \in \integer^+$).

 \end{itemize}

Without loss of generality:

\begin{enumerate} \renewcommand{\theenumi}{\alph{enumi}}

\item The action is orientation preserving.

\item We may assume that $U_1 = U$ (because $U_1$ is contained in a maximal unipotent subgroup of $\SL(2,\rational)$, and all maximal unipotent subgroups of $\SL(2,\rational)$ are conjugate).

\item  \label{SL(2,Z[1/r])UBddOrbPf-Uunbdd}
 We may assume that the $U$-orbit of~$x$ is \textbf{not} bounded above. (Otherwise, it would not be bounded below, and we could reverse the orientation of~$\real$.)

\item  \label{SL(2,Z[1/r])UBddOrbPf-rank1}
We may assume 
$$\lim_{\begin{matrix}
 u \to +\infty \\  \umat{u} \in U
 \end{matrix}}
  x \cdot \umat{u} = \infty .$$
  (See~\ref{justify}(\theenumi).)

\item  \label{SL(2,Z[1/r])UBddOrbPf-Vunbdd}
We may assume 
$$\lim_{\begin{matrix}
 v \to +\infty \\  \vmat{v} \in V
 \end{matrix}}
  x \cdot \vmat{v} = \infty .$$
  (This would be obvious from~\pref{SL(2,Z[1/r])UBddOrbPf-rank1} if $V$ were conjugate to~$U$ in~$\Gamma$. In the general case, a bit of work is required (see~\ref{justify}(\theenumi)).)

\item \label{SL(2,Z[1/r])UBddOrbPf-fixedpt}
 We may assume that $\dmat{\omega}$ fixes~$x$. (It is not difficult to see that $\dmat{\omega}$ has a fixed point in the interval $[x,\infty)$ (see~\ref{justify}(\theenumi)), and there is no harm in replacing $x$ with this fixed point.)

\end{enumerate}
From \pref{SL(2,Z[1/r])UBddOrbPf-rank1}, we know there is some $u \in \integer[1/r]^+$, such that $x \cdot \vmat{1} < x \cdot \umat{u}$. Then, because the action of~$\Gamma$ is orientation preserving, we have
 $$ x \cdot \vmat{1} \dmat{\omega}^n < x \cdot \umat{u} \dmat{\omega}^n $$
for all $n \in \integer$.
On the other hand, as $n \to \infty$, we have
 $$x \cdot \vmat{1} \dmat{\omega}^n
 = (x \cdot \dmat{\omega}^n) \cdot (\dmat{\omega}^{-n} \vmat{1} \dmat{\omega}^n)
 = x \cdot \vmat{\omega^{2n}}
 \to +\infty
, $$
and
$$ \text{$x \cdot \umat{u} \dmat{\omega}^n
= (x \cdot \umat{\omega}^n) \cdot (\dmat{\omega}^{-n} \umat{u} \dmat{\omega}^n)
= x \cdot \umat{\omega^{-2n} u}
< x \cdot \umat{u}
$ is bounded} .$$
 This is a contradiction.
 \end{proof}
 
 Because B.~Liehl \cite{LiehlBddGenSL2} proved that $\SL \bigl( 2, \integer[1/r] \bigr)$ is boundedly generated by unipotents (or see \cite{CarterKellerPaige} or \cite{MorrisBddGenSL2}), the above proposition has the following consequence:

\begin{cor} \label{SL(2,Z[1/r])NotRO}
 Let\/ $\Gamma$ be a finite-index subgroup of\/ $\SL \bigl( 2, \integer[1/r] \bigr)$, for some natural number~$r>1$. Then:
 \begin{enumerate}
 \item $\Gamma$ has no nontrivial, orientation-preserving action on\/~$\real$,
 and
 \item $\Gamma$ is not right orderable.
 \end{enumerate}
 \end{cor}

This immediately implies the following generalization. (The noncompactness of the semisimple groups was omitted from the hypotheses by mistake in \cite{LifschitzMorris-CR}.)

\begin{cor}[{\cite[Thm.~1.1]{LifschitzMorris-CR}}] \label{SarithThm}
Suppose 
 \begin{itemize}
 \item $G_\infty$ is a connected, noncompact, real, semisimple Lie group with finite center, 
 \item $S$ is a finite, nonempty set of prime numbers,
 \item $G_p$ is a Zariski-connected, noncompact, semisimple algebraic group over the $p$-adic field\/~$\rational_p$, for each $p \in S$, 
 \item $G$ is isogenous to\/ $\bigtimes_{p \in S \cup \{\infty\}} G_p$,
 and
 \item $\Gamma$ is a noncocompact, irreducible lattice in~$G$.
 \end{itemize}
Then:
 \begin{enumerate}
 \item $\Gamma$ has no nontrivial, orientation-preserving action on\/~$\real$,
 and
 \item $\Gamma$ is not right orderable.
 \end{enumerate}
\end{cor}

\begin{proof}
The Margulis Arithmeticity Theorem \cite[Thm.~A, p.~298]{MargulisBook} tells us that $\Gamma$ must be an $S$-arithmetic subgroup of~$G$. This means there is an algebraic number field~$\F$, a semisimple algebraic group~$\G$ over~$\F$, and a finite set~$\scrS$ of places of~$\F$, such that (a finite-index subgroup of)~$\Gamma$ is isomorphic to a finite-index subgroup of $\G(\ints_{\scrS})$ (where $\ints_{\scrS}$ is the ring of $\scrS$-integers of~$\F$.
 \begin{itemize}
 \item Since $\Gamma$ is noncocompact, we must have $\Frank \G \ge 1$, so $\G$ contains a subgroup that is isogenous to $\SL(2, {\cdot})$.
 \item Since $\scrS$ is nonempty, there is a (rational) prime~$p$, such that $\integer[1/p] \subseteq \ints_{\scrS}$. 
 \end{itemize}
 Therefore, $\Gamma$ contains a subgroup that is commensurable to $\SL \bigl( 2, \integer[1/p] \bigr)$. So Proposition~\ref{SL(2,Z[1/r])UBddOrb} applies.
 \end{proof}
 
The following conclusion is obtained by combining the above corollary with a generalization of Ghys' Theorem \cite{GhysAOC} to the setting of $S$-arithmetic groups \cite[Cor.~6.11]{WitteZimmer-CircleBundle}:

\begin{cor} \label{SarithNoActS1}
Assume
\begin{itemize}
\item the hypotheses of Corollary~\ref{SarithThm},
and
\item no simple factor of~$G_\infty$ is isogenous to\/ $\SL(2,\real)$.
\end{itemize}
Then any action of\/~$\Gamma$ on the circle~$S^1$ factors through a finite quotient of\/~$\Gamma$.
\end{cor}

\begin{noheading}[\bf Justification of the assumptions in the proof of Proposition~\ref{SL(2,Z[1/r])UBddOrb}] \label{justify} \ 
 \begin{enumerate}
 
 \item[\pref{SL(2,Z[1/r])UBddOrbPf-rank1}] \label{justify-rank1}
 It is easy to see that the additive group of $\integer[1/r]$ has only two total orderings (such that $u_1 \prec u_2 \Rightarrow u_1 + u_3 \prec u_2 + u_3$); namely,
 $$ \text{either
  \quad
   $u_1 \prec u_2 \Leftrightarrow u_1 < u_2$
 \quad or \quad
  $u_1 \prec u_2 \Leftrightarrow -u_1 < -u_2$.} $$
   This implies that
 either
 \begin{enumerate} \renewcommand{\theenumii}{\roman{enumii}}
 \item  \label{justify-rank1-goodorder}
  $x \cdot \umat{u_1} < x \cdot \umat{u_2} \Leftrightarrow u_1 < u_2$,
 or
  \item  $x \cdot \umat{u_1} < x \cdot \umat{u_2} \Leftrightarrow -u_1 < -u_2$.
  \end{enumerate}
  (This conclusion can also be obtained as a special case of Corollary~\fullref{realhomo}{exists} below.)
  We may assume \pref{justify-rank1-goodorder} holds (by conjugating by $\begin{bmatrix} -1 & 0 \\ 0 & 1 \end{bmatrix}$ if necessary).
 The desired conclusion now follows from~\pref{SL(2,Z[1/r])UBddOrbPf-Uunbdd}.
 
  \item[\pref{SL(2,Z[1/r])UBddOrbPf-Vunbdd}] \label{justify-Vunbdd}
  There is some element~$\gamma$ of~$\Gamma$ that does not normalize~$U$ (since $U$ is not normal in~$\Gamma$). Then $\gamma$ does not normalize~$\U$, so $\U \neq \gamma^{-1} \U \gamma$. Since $\Qrank \bigl( \SL(2,\cdot) \bigr) = 1$, this implies some element~$g$ of $\SL(2,\rational)$ conjugates the pair $(\U,\gamma^{-1} \U \gamma)$ to the pair $(\U,\V)$ \see{OppUnipsConjInRank1}. Thus, replacing $\Gamma$ with $(g^{-1} \Gamma g) \cap \SL \bigl( 2, \integer[1/r] \bigr)$, and letting $\widehat\Gamma = g^{-1} \Gamma g$, we may assume $\Gamma$ is contained in a subgroup $\widehat\Gamma$ of $\SL(2,\rational)$, such that
   \begin{enumerate} \renewcommand{\theenumii}{\roman{enumii}}
   \item the action of~$\Gamma$ on~$\real$ extends to an orientation-preserving action of $\widehat\Gamma$ on~$\real$,
   and
   \item $g^{-1} \widehat U g = \widehat V$ for some $g \in \widehat\Gamma$ (where $\widehat U = \U \cap \widehat\Gamma$ and $\widehat V = \V \cap \widehat\Gamma$).
   \end{enumerate}
By generalizing \pref{SL(2,Z[1/r])UBddOrbPf-rank1} to the subgroup~$\widehat U$ and noting that $g$ preserves orientation, we see that 
$$ \lim_{\begin{matrix}
 u \to +\infty \\  \umat{u} \in \widehat U
 \end{matrix}}
  xg \cdot (g^{-1} \umat{u} g) = \infty .$$
 Because
 \begin{equation} \label{upperlowerconjugate}
  \begin{bmatrix} 0 & 1 \\ -1 & 0 \end{bmatrix}^{-1} \umat{u} \begin{bmatrix} 0 & 1 \\ -1 & 0 \end{bmatrix}
 =  \begin{bmatrix} 1 & 0 \\ -u & 1 \end{bmatrix}
 = \vmat{u} ,
 \end{equation}
 we see that this implies 
$$ \lim_{\begin{matrix}
 v \to \infty \\  \vmat{v} \in \widehat V
 \end{matrix}}
  xg \cdot \vmat{v} = \infty .$$
  Replacing $x$ with $\max\{x,x g\}$ (and restricting to the subgroup~$V$ of~$\widehat V$) yields the desired conclusion.

 \item[\pref{SL(2,Z[1/r])UBddOrbPf-fixedpt}] \label{justify-fixedpt}
 We wish to show that $\dmat{\omega}$ has a fixed point in the interval $[x,\infty)$; that is, we wish to show that the orbit of~$x$ under the group $\langle \dmat{\omega} \rangle$ is bounded above. For definiteness, let us assume that $x \cdot \dmat{\omega} \ge x$. (This causes no loss of generality, because the transpose-inverse automorphism of $\SL \bigl( 2, \integer[1/r])$ sends $\dmat{\omega}$ to its inverse while interchanging~$U$ with~$V$.) From \pref{SL(2,Z[1/r])UBddOrbPf-rank1}, we know there is some $u > 0$, such that $x \cdot \dmat{\omega} < x \cdot \umat{u}$. For convenience, let 
 $$u_n = u(1 + \omega^{-2} + \cdots + \omega^{-2n}) \in \integer[1/r] ,$$
 so
 $$ u_n = u + \omega^{-2} u_{n-1} .$$
  Then, by induction on~$n$, we have
  \begin{align*}
   x \cdot \dmat{\omega}^n
  &= (x \cdot \dmat{\omega}^{n-1}) \cdot \dmat{\omega}
  < ( x \cdot \umat{u_{n-1}}) \cdot \dmat{\omega}
  \\&= ( x \cdot \dmat{\omega})  \cdot ( \dmat{\omega}^{-1} \umat{u_{n-1}} \dmat{\omega} )
  < ( x \cdot \umat{u}) \cdot ( \dmat{\omega}^{-1} \umat{u_{n-1}} \dmat{\omega} )
 = x \cdot \umat{u_{n}} 
 . \end{align*}
 Since the geometric series $\{u_n\}$ converges (hence is bounded above), we conclude that $\{ x \cdot \dmat{\omega}^n\}$ is bounded above.
 
 \end{enumerate}

\end{noheading}

\begin{rem}
The main difficulty in proving the result with a ring~$\ints$ of integers in the place of $\integer[1/r]$ is that the additive group of~$\ints$ has infinitely many different orderings. (It is isomorphic to $\integer^k$, for some $k > 1$, and any faithful homomorphism to~$\real$ yields an ordering.) Because of this, the natural analogue of assumption \pref{SL(2,Z[1/r])UBddOrbPf-rank1} is not at all obvious. By using the fact that $U$ is normalized by~$\dmat{\omega}$, it will be shown that only finitely many orderings of~$\ints$ can arise (see Step~\ref{pEigenvalue} on page~\pageref{pEigenvalue}). It is then easy to adapt the proof of Proposition~\ref{SL(2,Z[1/r])UBddOrb} to apply to~$\SL(2,\ints)$.
\end{rem}

\section{Preliminaries on arithmetic groups} \label{PrelimArithSect}

We recall some well-known facts.

\begin{thm}[{(Margulis Normal Subgroup Theorem \cite[(A), p.~258]{MargulisBook})}]
\label{MargNormSubgrp}
 If
 \begin{itemize}
 \item $G$ is a connected, semisimple Lie group with trivial center,
 \item $\Rrank G \ge 2$,
 \item $\Gamma$ is an irreducible lattice in~$G$,
 and
 \item $N$ is a nontrivial, normal subgroup of\/~$\Gamma$, 
 \end{itemize}
 then\/ $\Gamma/N$ is finite.
 \end{thm}
 
 \begin{cor} \label{ActMustBeFaithful}
  If\/ $\Gamma$ is as in Theorem~\ref{MargNormSubgrp}, then any nontrivial, orientation-preserving action of\/~$\Gamma$ on\/~$\real$ is faithful.
   \end{cor}

\begin{proof}
If an action of~$\Gamma$ is not faithful, then its kernel is a nontrivial, normal subgroup of~$\Gamma$. Hence, the kernel has finite index, so $\Gamma$ acts via homeomorphisms of finite order. Since $\real$ has no nontrivial, orientation-preserving homeomorphisms of finite order, we conclude that the action is trivial.
\end{proof}

 \begin{cor} \label{Act<>RO}
  Suppose\/ $\Gamma$ is as in Theorem~\ref{MargNormSubgrp}. Then\/ $\Gamma$ is right orderable if and only if\/ $\Gamma$ has a nontrivial, orientation-preserving action on\/~$\real$.
   \end{cor}

\begin{proof} 
 It is well known that a countable group is right orderable if and only if it has a faithful, orientation-preserving action on~$\real$ \cite[Thm.~6.8]{GhysSurvey}. 
\end{proof}

\begin{prop} \label{NonsimpleContainsSL2}
Suppose 
 \begin{itemize}
 \item $G$ is a connected, noncompact, semisimple Lie group with finite center, 
 \item $\Rrank G \ge 2$,
 \item the adjoint group of~$G$ is \textbf{not} simple,
 and
 \item $\Gamma$ is a noncocompact, torsion-free, irreducible lattice in~$G$.
 \end{itemize}
Then some subgroup of\/~$\Gamma$ is isomorphic to a finite-index subgroup of\/ $\SL(2,\ints)$, where $\ints$ is the ring of integers of a number field\/~$\F$ that is neither\/~$\rational$ nor an imaginary quadratic extension of\/~$\rational$.
\end{prop}

\begin{proof}
The Margulis Arithmeticity Theorem \cite[Thm.~A, p.~298]{MargulisBook} tells us that $\Gamma$ must be an arithmetic subgroup of~$G$. This means there is an absolutely almost-simple algebraic group~$\G$ over some algebraic number field~$\F$, such that 
\begin{itemize}
\item a finite-index subgroup of~$\Gamma$ is isomorphic to a finite-index subgroup of $\G(\ints)$ (where $\ints$ is the ring of integers of~$\F$),
and
\item $G$ is isogenous to $\bigtimes_{v \in S_\infty} G(\F_v)$, where $S_{\infty}$ is the set of infinite places of~$\F$.
\end{itemize}
(The assumption that $\Gamma$ is noncocompact implies that each $G(\F_v)$ is noncompact.)
Since $\Gamma$ is noncocompact, we must have $\Frank \G \ge 1$, so $\G$ contains a subgroup that is isogenous to $\SL(2, {\cdot})$.  Therefore, $\Gamma$ contains a subgroup that is commensurable to $\SL( 2, \ints)$. Since the adjoint group of~$G$ is not simple, we know that the product $\bigtimes\nolimits_{v \in S_\infty} G(\F_v)$ has more than one factor, so $\F$~has more than one infinite place. Therefore, $\F$ is neither~$\rational$ nor an imaginary quadratic extension of~$\rational$.
  \end{proof}

The following observation is a simple case of much more general superrigidity
theorems \cite{GorbacevicSuper, StarkovSuper, WitteSuperLatt}. Its proof can
be reduced to the abelian case by using the fact that $[\Gamma,\Gamma]$ is a
lattice in $[\U,\U]$.

\begin{lem}[{(cf.\ \cite[Thm.~2.11, p.~33]{Raghunathan})}] \label{MalcevRig}
 If
 \begin{itemize}
 \item $\U$ is a $1$-connected, nilpotent Lie group,
 \item $\Gamma$ is a lattice in~$\U$,
 and
 \item $\sigma \colon \Gamma \to \real$ is any homomorphism,
 \end{itemize}
 then $\sigma$ extends uniquely to a continuous homomorphism $\hat\sigma
\colon \U \to \real$.
 \end{lem}
  
 The following well-known property of groups of rank~1 is useful.
 
 \begin{prop}[{(cf.\ \cite[(8.4) and (4.8), pp.~124 and~88]{BorelTits})}] \label{OppUnipsConjInRank1}
 Suppose
 \begin{enumerate}
 \item $\G$ is a semisimple algebraic $\rational$-group,
\item $\Qrank G = 1$,
\item $\GU_1$, $\GV_1$, $\GU_2$, and\/ $\GV_2$ are maximal unipotent subgroups of\/~$\G$,
and
\item $\GU_i \neq \GV_i$ for $i = 1,2$.
\end{enumerate}
Then there exists $g \in \G(\rational)$, such that $g^{-1} \GU_1 g = \GU_2$ and $g^{-1} \GV_1 g = \GV_2$.
\end{prop}

\section{Preliminaries on bounded generation} \label{PrelimBddGenSect}

The following observation shows that being virtually boundedly generated by unipotents is not affected by passing to a finite-index subgroup.

\begin{lem}[{(cf.\ \cite[Prop.\ on p.~256]{Murty-BddGen})}]
\label{FinIndHasBddGen}
 Suppose
 \begin{itemize}
 \item $\Gamma$ is a group,
 \item $U_1,\ldots,U_m$ are subgroups of\/~$\Gamma$,
 \item the product $U_1 U_2 \cdots U_m$ is a finite-index subgroup of\/~$\Gamma$,
 and
 \item $\Gamma'$ is a finite-index subgroup of\/~$\Gamma$.
 \end{itemize}
 Then there is a list  $U_1',\ldots,U_n'$ of finitely many subgroups of\/~$\Gamma'$, such that
 \begin{enumerate}
 \item each $U'_i$ is conjugate {\rm(}in\/ $\Gamma${\rm)} to a subgroup of some $U_{j_i}$,
 and
  \item the product $U'_1 U'_2 \cdots U'_m$ is a finite-index subgroup of\/~$\Gamma'$.
  \end{enumerate}
  \end{lem}

\begin{proofnum}{Lemma~\ref{BddGen+BddOrbits}} \label{BddGen+BddOrbitsPf}
 Suppose we are given a nontrivial, orientation-preserving action of~$\Gamma$
on~$\real$. (This will lead to a contradiction.) Because the action is
nontrivial, some point~$x_0$ of~$\real$ is \emph{not} fixed by all
of~$\Gamma$. Let
 $$ \mbox{$a = \inf (x_0 \cdot \Gamma)$ \qquad and \qquad $b = \sup (x_0 \cdot \Gamma)$,} $$
 where $x_0 \cdot \Gamma$ denotes the $\Gamma$-orbit of~$x_0$.
 Then $(a,b)$ is a (nonempty) $\Gamma$-invariant open interval, so it is homeomorphic to~$\real$. By replacing
$\real$ with this subinterval, we may assume that
 $$ \mbox{$\inf (x_0 \cdot \Gamma) = - \infty$ \qquad and \qquad $\sup (x_0 \cdot
\Gamma) = \infty$;} $$
 thus, the $\Gamma$-orbit of~$x_0$ is \emph{not} bounded.

By passing to a finite-index subgroup, we may assume that 
 $$ \Gamma = U_1 U_2 \cdots U_m .$$
By induction on~$m$, we see
that $x_0 \cdot \Gamma = x_0 \cdot (U_1 U_2 \cdots U_m)$ is bounded. This contradicts the
conclusion of the preceding paragraph.
 \end{proofnum}

\section{Preliminaries on unbounded orbits of unipotent subgroups} \label{PrelimUnbddOrbitSect}

Let us recall the following fundamental result on right-orderings of
nilpotent groups that was proved (independently) by J.~C.~Ault \cite{Ault}
and A.~H.~Rhemtulla \cite{Rhemtulla}. We state only a weak version.

\begin{thm}[{(Ault, Rhemtulla \cite[Thm.~7.5.1, p.~141]{MuraRhemtulla})}]
\label{AultRhemtulla}
 Suppose 
 \begin{itemize}
 \item $U$ is a finitely generated, nilpotent group,
 and
 \item we have an orientation-preserving action of~$U$ on\/~$\real$, such that
$0$ is \textbf{not} a fixed point.
 \end{itemize}
 Then there is a nontrivial homomorphism $p \colon U \to \real$, such that
 \begin{equation}
 \mbox{for all $u \in U$ with $p(u) > 0$, we have\/  $0 \cdot u > 0$.}
 \end{equation}
 Furthermore,
\begin{enumerate}
\item $p$~is unique, up to multiplication by a positive scalar.
 \item \label{AultRhemtulla-p=0}
 If $z,u \in U$ with $p(z) = 0 < p(u)$, then $z$ has a fixed point in the closed interval\/ $[0, 0 \cdot u]$.
 \end{enumerate}
 \end{thm}

\begin{rem} \label{AultRhemtulla>}
For all $u_1,u_2 \in U$ with $p(u_1) < p(u_2)$, we have  $0 \cdot u_1 < 0 \cdot u_2$.
This is because $p(u_2 u_1^{-1}) = p(u_2) - p(u_1) > 0$, so $0 \cdot u_2 u_1^{-1} > 0$.
\end{rem}

\begin{cor} \label{realhomo}
 Suppose 
 \begin{itemize}
 \item $\U$ is a $1$-connected, nilpotent Lie group,
 \item $U$ is a lattice in~$\U$,
 \item we have an orientation-preserving action of~$U$ on\/~$\real$, 
 and
 \item the $U$-orbit of\/~$0$ is \textbf{not} bounded above.
 \end{itemize}
 Then:
 \begin{enumerate}
 \item \label{realhomo-exists}
 There is a nontrivial, continuous homomorphism $p_{U} \colon \U \to
\real$, such that, for 
\begin{itemize}
 \item all $u_1,u_2 \in U$ with $p_{U}(u_1) <
p_{U}(u_2)$, 
and
\item all $x \ge 0$,
\end{itemize}
we have  $x \cdot u_1 < x \cdot u_2$. 
 \item \label{realhomo-unique}
 The homomorphism~$p_{U}$ is unique up to multiplication by a positive
scalar.
 \end{enumerate}
 \end{cor}

\begin{proof}
\pref{realhomo-unique} The uniqueness of $p_{U}$ is a consequence of the
uniqueness in \pref{AultRhemtulla} and \pref{MalcevRig}.

 \pref{realhomo-exists} Let
 \begin{itemize}
\item $p \colon U \to \real$ be the homomorphism provided by the
Ault-Rhemtulla Theorem \pref{AultRhemtulla}. 
\item $p_{U} \colon \U \to \real$ be the (unique) continuous homomorphism that extends~$p$  \see{MalcevRig}.
\item $C$ be the component of
 $$ \{\, x \in \real \mid \mbox{for all $u \in U$ with $p_U(u) > 0$, we have $x \cdot u > x$} \,\} ,$$
that contains~$0$.
 \item $y$ be the upper endpoint of the interval~$C$.
 \end{itemize}
 It suffices to show that $y = \infty$ \cf{AultRhemtulla>}.

Let us suppose $y$ is finite. (This will lead to a contradiction.)
 Since the $U$-orbit of~$0$ is not bounded above, we know that $U$ has no fixed points in $[0,\infty)$; therefore, $y$ is not a fixed point. Thus, combining the Ault-Rhemtulla Theorem (with $y$ in the role of~$0$) with \pref{MalcevRig} yields a nontrivial, continuous homomorphism $p'_U \colon \U \to \real$, such that
 for all $u \in U$ with $p'_U(u) > 0$, we have  $y \cdot u > y$.
Because $U$ acts continuously on~$\real$, we have 
 $$ \{\, u \in U \mid p_{U}(u) > 0 \,\}
  \subseteq \{\, u \in U \mid y \cdot u \ge y \,\}
  \subseteq \{\, u \in U \mid p'_{U}(u) \ge 0 \,\}
   .$$
 This implies that $p'_{U}(u) = p_{U}(u)$ (up to a positive scalar
multiple). 
 \begin{itemize}
 \item Fix some $u_0 \in U$ with $p_{U}(u_0) > 0$. 
 \item Consider any $x \in [y, y \cdot u_0)$.
 \end{itemize}
 For any $u \in U$ with $p_{U}(u) > 0$, there is some positive
integer~$k$, such that $p_{U}(u^k) > p_{U}(u_0)$. Thus, 
 $$ x \cdot u^k \ge y \cdot u^k > y \cdot u_0 > x ,$$
 so $x \cdot u > x$. Since $x$ is an arbitrary element of $[y,y \cdot u_0)$, we conclude that $[y,y \cdot u_0) \subseteq C$. 
 
 This contradicts the fact that $y$ is the upper endpoint of~$C$. 
 \end{proof}
 
 \begin{notation}
 Suppose we are given an orientation-preserving action of a group~$\Gamma$ on~$\real$.
 For convenience in the remaining proofs of this section, we define a partial order~$\prec$ on~$\Gamma$ by
  $$ \mbox{$g \prec h$ \qquad $\Leftrightarrow$ \qquad $0 \cdot g \< 0 \cdot h$.} $$
 \end{notation}

\begin{cor} \label{p=root}
 Suppose 
 \begin{itemize}
 \item $\G$ is an almost simple algebraic\/ $\rational$-group,
 \item $\Gamma$ is an arithmetic subgroup of\/~$\G$,
 \item we are given an orientation-preserving action of\/~$\Gamma$ on\/~$\real$,
 \item $\GU$ is a unipotent\/ $\rational$-subgroup
of~$\G$,
 \item $\GT$ is a\/ $\rational$-torus of\/~$\G$ that normalizes\/ $\GU$,
  \item $U = \GU \cap \Gamma$ and $T = \GT \cap \Gamma$,
 \item the $U$-orbit of\/~$0$ is not bounded above,
 and
 \item $p_U \colon \GU \to \real$ is as specified in Corollary~\fullref{realhomo}{exists}.
 \end{itemize}
 Then there exist
 \begin{itemize}
 \item a one-parameter\/ $\real$-subgroup\/ $\GU_1$ of\/~$\GU$ that is normalized
by~$T$,
 and
 \item a $T$-equivariant projection $\pi \colon \GU_\real \to \GU_1$, 
 \end{itemize}
 such that 
 $$ p_U(u) = p_U \bigl( \pi(u) \bigr) $$
 for all $u \in \GU_{\real}$.
 \end{cor}

\begin{proof}
 It suffices to show that $T$ normalizes the kernel of~$p_U$. (Because $T$ is
a torus, this implies there is a complementary subgroup~$\GU_1$ that is
normalized by~$T$.) 

Suppose some $t \in T$ does not normalize the kernel of~$p_U$. (This will
lead to a contradiction.) 
We may assume, without loss of generality, that $e \prec t$ (by replacing $t$ with~$t^{-1}$, if necessary).
It is not difficult to see there must be some $u \in U$, such that $p_U(u) > 0$, but
 \begin{equation} \label{p=rootPf-p<0}
 p_U(t u t^{-1}) < 0 .
 \end{equation}
 There is some $v \in U$, such that $t \prec v$ (because the $U$-orbit of~$0$ is not bounded above). Choose a large integer $k>0$, so
that $k \, p_U(u) > p_U(v)$; thus, $e \prec u^k v^{-1}$. Then
 $$ e \prec t \, (u^k v^{-1}) \, (v t^{-1}) = t u^k t^{-1} = (t u t^{-1})^k
,$$
 so $e \prec t u t^{-1}$. This contradicts \pref{p=rootPf-p<0}.
 \end{proof}

\begin{cor} \label{OppScalarsSame}
 Suppose 
 \begin{itemize}
 \item $\G$ is an almost simple algebraic\/ $\rational$-group,
 \item $\Gamma$ is an arithmetic subgroup of\/~$\G$,
 \item we are given an orientation-preserving action of\/~$\Gamma$ on\/~$\real$,
 \item $\GU$ and\/ $\GV$ are unipotent\/ $\rational$-subgroups
of\/~$\G$,
 \item $\GT$ is a\/ $\rational$-torus of\/~$\G$ that normalizes both\/ $\GU$ and\/~$\GV$,
  \item $U = \GU \cap \Gamma$, $V = \GV \cap \Gamma$, and $t \in \GT \cap
\Gamma$,
 and
 \item the $U$-orbit of\/~$0$ and the $V$-orbit of\/~$0$ are not bounded above.
 \end{itemize}
 Then:
 \begin{enumerate}
 \item \label{OppScalarsSame-scalars}
 There are real scalars $\omega_U$ and~$\omega_V$, such that 
 $$ \mbox{$p_U(t^{-1} u t) = \omega_U \, p_U(u)$\quad and\quad $p_V(t^{-1} v
t) = \omega_V \, p_V(v)$,} $$
 for all $u \in U$ and $v \in V$.
 \item \label{OppScalarsSame-FP}
 If\/ $|\omega_U| \neq 1$ and\/ $|\omega_V| \neq 1$, then $t$ fixes some
point~$x$ of\/~$\real$, such that the $U$-orbit of~$x$ and the $V$-orbit of~$x$
are not bounded above.
 \item \label{OppScalarsSame-same}
 If\/ $|\omega_U| < 1$, then\/ $|\omega_V| \le 1$.
 \end{enumerate}
 \end{cor}

\begin{proof}
 \pref{OppScalarsSame-scalars}
 This follows from Corollary~\ref{p=root}.

\medskip

\pref{OppScalarsSame-FP}
 Let $F_U$ and $F_V$ be the fixed-point sets of~$U$ and~$V$, respectively,
and let 
 $$\alpha = \max \{ \sup F_U, \sup F_V \} .$$
 (Note that $\alpha < 0$, because the $U$-orbit of~$0$ and the $V$-orbit
of~$0$ are not bounded above.) It suffices to show that $t$ has a fixed point
in $(\alpha,\infty)$.

Assume, without loss of generality, that $\alpha = \sup F_U$. Because $t$
normalizes $U$, we know that $F_U$ is $t$-invariant, so the interval
$(\alpha,\infty)$ is $t$-invariant. By replacing $\real$ with this interval
(and ignoring~$V$), we may assume $\alpha = -\infty$. Thus, 
 \begin{itemize}
 \item the $U$-orbit of~$0$ is neither bounded below nor bounded above,
 and
 \item it suffices to show that $t$ has a fixed point (anywhere in~$\real$).
 \end{itemize}

 We may assume 
 \begin{itemize}
 \item $\omega_U > 0$, by replacing $t$ with~$t^2$,
 \item $\omega_U < 1$, by replacing $t$ with~$t^{-1}$, if necessary,
 and
 \item $t \succ e$, by reversing the orientation of~$\real$, if necessary.
 \end{itemize}
 Because the $U$-orbit of~$0$ is not bounded above, there is some $u \in U$
with 
$$t \prec u .$$
 Choose some $u_0 \in U$ with 
 $$p(u_0) > \frac{p(u)}{1 - \omega_U} .$$
 Then, for every $k > 0$, we have
 \begin{align*}
 &p_U \bigl( (t^{-0} u t^0) (t^{-1} u t^1) (t^{-{2}} u t^{2}) \cdots 
(t^{-(k-1)} u t^{k-1}) \bigr) \\
 &\hskip0.5in =   p_U(u) \, (1 + \omega_U + \omega_U^2 + \cdots + \omega_U^{k-1}) \\
 &\hskip0.5in < p(u_0) ,
 \end{align*}
 so
 $$ (t^{-0} u t^0) (t^{-1} u t^1) (t^{-{2}} u t^{2}) \cdots (t^{-(k-1)} u
t^{k-1})
 \prec u_0 .$$ 
 Therefore
 \begin{align*}
  e &\prec (u t^{-1})^k \\
 &= u \, (t^0 t^{-1}) \, u \, (t^1 t^{-2}) \, u \, (t^2 t^{-3})
\cdots u \, (t^{k-2} t^{-(k-1)}) \, u \, (t^{k-1} t^{-k}) \\
 &= (t^{-0} u t^0) (t^{-1} u t^1) (t^{-{2}} u t^{2}) \cdots  (t^{-(k-1)} u
t^{k-1}) t^{-k} \\
 &\prec u_0 t^{-k}
 .
 \end{align*}
 This means
 $ t^k \prec u_0$ (for every~$k$),
 so the $\langle t \rangle$-orbit of~$0$ is bounded above (by~$0 \cdot u_0$).
Therefore $t$ has a fixed point.

\medskip

\pref{OppScalarsSame-same}
 Suppose $|\omega_U| \le 1$ and $|\omega_V| > 1$. (This will lead to a
contradiction.) By replacing $t$ with~$t^2$, we may assume $\omega_U$
and~$\omega_V$ are positive. We may also assume that $t$~fixes~$0$
(see~\pref{OppScalarsSame-FP}). Fix 
 \begin{itemize}
 \item $v \in V$ with $p_V(v) > 0$,
 and
 \item $u \in U$ with $0 \cdot v \< 0 \cdot u$ and $p_U(u) > 0$.
 \end{itemize}
 Because the action of~$\Gamma$ is orientation preserving, we have
 \begin{equation} \label{OppScalarsSame-same-vtn<utn}
  0 \cdot v t^n < 0 \cdot u t^n  
  \end{equation}
for all $n \in \integer$.
Note that, as $n \to \infty$, we have $p_V(t^{-n} v t^n) = \omega_V^n p_V(v) \to +\infty$, so
 $$0 \cdot v t^n
 = (0 \cdot t^n) \cdot (t^{-n} v t^n)
 = 0 \cdot (t^{-n} v t^n)
 \to +\infty
.$$
On the other hand, we have $p_U(t^{-n} u t^n) = \omega_U^n p_U(u) \to 0 < p_U(u)$, so
$$ \text{$0 \cdot u t^n
= (0 \cdot t^n) \cdot (t^{-n} u t^n)
= 0 \cdot (t^{-n} u t^n)
< 0 \cdot u$ is bounded} .$$
 This contradicts \pref{OppScalarsSame-same-vtn<utn}.
 \end{proof}

\section{Proof of Theorem~\ref{ArchimedeanUBddOrbits}} \label{ArchimedeanUBddOrbitsPf}

 Throughout this section, the conditions in the statement of Theorem~\ref{ArchimedeanUBddOrbits} are satisfied:
\begin{itemize}
\item $\F$ is an algebraic number field that is neither~$\rational$ nor an imaginary quadratic extension of~$\rational$,
\item $\ints$ is the ring of integers of~$\F$,
\item no proper subfield of~$\F$ contains a finite-index subgroup of~$\ints^\times$,
and
\item $\Gamma$ is a finite-index subgroup of $\SL(2,\ints)$.
\end{itemize}
 Furthermore, we are given an orientation-preserving action of~$\Gamma$ on~$\real$. We wish to show that every orbit of every unipotent subgroup of~$\Gamma$ is bounded.

\begin{notation} \ 

 \begin{itemize}
  
 \item For $u,v,w \in \complex$, with $w \neq 0$, let 
 $$ \mbox{$\umat{u} = \begin{bmatrix}
 1 & u \\
 0 & 1 \\
 \end{bmatrix}$,
 \qquad
 $\vmat{v} = \begin{bmatrix}
 1 & 0 \\
 -v & 1 \\
 \end{bmatrix}$
  \quad and \quad
 $\dmat{w} = \begin{bmatrix} w & 0 \\ 0 & 1/w \end{bmatrix}$%
 .}$$
 
 \item Let
 $$ \mbox{$\U = \bigset{ \umat{u} }{ u \in \F }$
 \quad and \quad
 $\V = \bigset{ \vmat{v} }{ v \in \F }$},$$
 so $\U$ and~$\V$ are opposite maximal unipotent subgroups of
$\SL(2,\F)$. 

\item Let
 $$ \mbox{$U = \U \cap \Gamma$ \quad and \quad $V = \V \cap \Gamma$.}$$

 \end{itemize}
 
 \end{notation}

\begin{assump}
 Assume some orbit of some unipotent subgroup~$U_0$ of~$\Gamma$ is not bounded.
(This will lead to a contradiction.) There is no harm in assuming, for definiteness, that:
\begin{enumerate}
\item $U_0 = U$ is a maximal unipotent subgroup,
 and 
 \item the $U$-orbit of~$0$ is not bounded above.
\end{enumerate}
\end{assump}

The proof now proceeds in a sequence of steps.

\setcounter{step}{0}

\begin{step} \label{ManyConjs}
 There is a sequence $g_1,g_2,g_3,\ldots$ of elements of\/~$\Gamma$, such
that
 \begin{enumerate} \renewcommand{\theenumi}{\alph{enumi}}
 \item \label{ManyConjs-distinct}
 the conjugates $g_1^{-1} U g_1, g_2^{-1} U g_2, g_3^{-1} U g_3,
\ldots$ are distinct, 
 and
 \item \label{ManyConjs-unbounded}
 for each~$j$, the $g_j^{-1} U g_j$-orbit of~$0$ is \textbf{not} bounded
above.
 \end{enumerate}
 \end{step}
 Because the normalizer $N_\Gamma(U)$ has infinite index in~$\Gamma$, there
is a sequence $g_1,g_2,g_3,\ldots$ of elements of~$\Gamma$, such that the
cosets $N_\Gamma(U) g_j$ are distinct. By passing to a subsequence (and
taking inverse of every term in the sequence, if necessary), we may assume
$0 \cdot g_j^{-1} \ge 0$, for each~$j$. 
\begin{enumerate}
\item[\pref{ManyConjs-distinct}]
The conjugates $g_1^{-1} U g_1, g_2^{-1} U g_2, g_3^{-1} U g_3, \ldots$
are distinct, because the cosets $N_\Gamma(U) g_j$ are distinct.

\item[\pref{ManyConjs-unbounded}]
Because $0 \cdot U$ is not bounded above, and $0 \cdot g_j^{-1} \ge 0$, it
is clear that $0 \cdot g_j^{-1} U$ is not bounded above. Therefore, $0
\cdot g_j^{-1} U g_j$ is not bounded above.
\end{enumerate}

\begin{step} \label{pEigenvalue}
 For each~$j$, define $p_{g_j^{-1} U g_j} \colon g_j^{-1} \U g_j \to \real$ as in Corollary~\fullref{realhomo}{exists}.
 There is a field embedding $\sigma_j \colon \F \hookrightarrow
\complex$, such that
  $$ p_{g_j^{-1} U g_j}\bigl( g_j^{-1}( \dmat{\omega}^{-1} u \dmat{\omega} ) g_j
\bigr) = \sigma_j(\omega)^{-2} \cdot p_{g_j^{-1} U g_j}( g_j^{-1} u g_j ) ,$$
 for every $u \in U$, and every $\omega \in \ints$, such
that $\dmat{\omega} \in \Gamma$.
 \end{step}
 To simplify the notation, assume, without loss of generality, that
$g_j = e$.
  Let 
 \begin{itemize}
 \item $S$ be the set of all archimedean places of~$\F$,
 \item $\F_\sigma$ be the completion of~$\sigma(\F)$, for each $\sigma \in
S$, so
 $$\F_\sigma = \begin{cases}
 \real & \mbox{if $\sigma(\F) \subset \real$} \\
 \complex & \mbox{if $\sigma(\F) \not\subset \real$}
 , \end{cases} $$
 \item $\G_S = \bigtimes_{\sigma \in S} \SL(2,\F_\sigma)$, 
 \item $\GU_S = \bigtimes_{\sigma \in S} \U_\sigma$, where $\U_\sigma = \{\, \umat{u} \mid u \in \F_\sigma\,\}$,
 \item $\GT_S = \bigtimes_{\sigma \in S} \T_\sigma$, where $\T_\sigma = \{\, \dmat{w} \mid w \in \F_\sigma, \ w \neq 0 \,\}$,
 and
 \item $\tuple{\phantom{x}} \colon \SL(2,\F) \hookrightarrow \G_S$ be defined
by 
 $\tuple g = \bigl( \sigma(g) \bigr)_{\sigma \in S}$.
 \end{itemize}
 It is well known from ``restriction of scalars" \cite[\S2.1.2, pp.~49--51]{PlatonovRapinchukBook} that $\G_S$ can be viewed as the $\real$-points of an almost simple algebraic $\rational$-group,
 such that
 \begin{itemize}
 \item $\tuple{\Gamma}$ is an arithmetic subgroup of~$\G_S$,
 and
 \item $\GU_S$ and $\GT_S$ are $\rational$-subgroups of~$\G_S$.
 \end{itemize}
 Let 
 \begin{itemize}
 \item $p_U \colon \GU_S \to \real$ be the homomorphism provided by Corollary~\ref{realhomo},
 \item $\GU_1$ be the one-parameter subgroup of~$\GU_S$ provided by Corollary~\ref{p=root},
 and
 \item $ \ints' = \{\, \omega \in \ints \mid \dmat{\omega} \in \Gamma\,\}$.
 \end{itemize}
 Since the Lie algebra of $\GU_1$ is a one-dimensional real subspace normalized by $\Ad \dmatCAP{\ints'}$, it must be contained in a single real eigenspace of $\Ad \dmat{\omega}$, for each $\omega \in \ints'$. Because no proper subfield of~$\F$ contains a finite-index subgroup of~$\ints^\times$, we know that
 \begin{itemize}
 \item no two distinct elements of~$S$ have the same restriction to~$(\ints')^2$, 
 and 
 \item $\sigma \bigl( (\ints')^2 \bigr) \not\subseteq \real$ whenever $\sigma$ is a complex place. 
\end{itemize}
Hence, $\GU_1$ must be contained in a single
factor~$\GU_{\sigma}$ of~$\GU_S$ (for some $\sigma \in S$), and $\sigma$ must be a real place. (Furthermore, the kernel of the projection~$\pi$ of Corollary~\ref{p=root} must contain $\GU_{\sigma'}$, for every $\sigma' \neq \sigma$.) Since $\GU_{\sigma} \iso \real$, there are only two nontrivial homomorphisms from $\GU_{\sigma}$ to~$\real$, up to multiplication by a positive scalar. Thus, we may assume
 $$ \mbox{$p_U \bigl( \umat{\alpha} \bigr) = \pm \sigma(\alpha)$ for $\alpha \in \F$}$$
 (and the same sign is used for all~$\alpha$).
 Because $\dmat{\omega}^{-1} \, \umat{\alpha} \, \dmat{\omega}
 = \umat{\omega^{-2} \alpha} $,
  we conclude that 
 $$ p_U \bigl( \dmat{\omega}^{-1} \, \umat{\alpha} \, \dmat{\omega} \bigr)
 = \pm\sigma \bigl({\omega^{-2} \alpha } \bigr)
 = \sigma(\omega)^{-2} \, p_U \bigl( \umat{\alpha} \bigr) ,$$
 as desired.

\begin{step} \label{ArchimedeanBddPf-pU=pVopp}
 We may assume
 \begin{enumerate} \renewcommand{\theenumi}{\alph{enumi}}
 \item \label{ArchimedeanBddPf-pU=pVopp-bdd}
 the $U$-orbit of\/~$0$ and the $V$-orbit of\/~$0$ are not bounded above,
 and
 \item \label{ArchimedeanBddPf-pU=pVopp-omega}
 we have
  $$ \mbox{$p_U ( \dmat{\omega}^{-1} u \dmat{\omega} ) = \omega^{-2} \cdot p_U(u)$
 \quad and \quad
 $p_V ( \dmat{\omega}^{-1} v \dmat{\omega}) = \omega^{2} \cdot p_V(v) $,}$$
 for all $u \in U$, all $v \in V$, and all $\omega \in \ints$, such that
$\dmat{\omega} \in \Gamma$.
 \end{enumerate}
 \end{step}
 Note that:
 \begin{itemize}
 \item Because there are only finitely many embeddings of~$\F$ in~$\complex$,
we may assume, by passing to a subsequence of~$\{g_j\}$, that $\sigma_1 =
\sigma_2$. 
 \item By replacing $\Gamma$ with $\sigma_1(\Gamma)$, we may assume
$\sigma_1$ is the natural inclusion $\sigma_1(\alpha) = \alpha$. 
 \item Because $\Frank \SL(2,\F) = 1$, we know that any pair of (unequal) maximal unipotent
$\rational$-subgroups of $\SL(2,\F)$ is conjugate to any other pair \see{OppUnipsConjInRank1}, so, by passing to a conjugate, we may assume $g_1^{-1} U g_1 = U$ and $g_2^{-1} U g_2 = V$ (cf.~\ref{justify}\pref{SL(2,Z[1/r])UBddOrbPf-Vunbdd}).
 \end{itemize}

\pref{ArchimedeanBddPf-pU=pVopp-bdd} From Step~\ref{ManyConjs}, we know that the
$U$-orbit of~$0$ and the $V$-orbit of~$0$ are not bounded above.
 
\pref{ArchimedeanBddPf-pU=pVopp-omega} The first half
of~\pref{ArchimedeanBddPf-pU=pVopp-omega} is immediate from
Step~\ref{pEigenvalue} (with $j = 1$).
 Taking 
 $$g_2 = \begin{bmatrix} 0 & 1 \\ -1 & 0 \end{bmatrix} ,$$
 noting that $g_2^{-1} \dmat{\omega} g_2 = \dmat{\omega}^{-1}$ and $\sigma_2 = \Id$,
 and letting $v = g_2^{-1} u g_2$, we see, from Step~\ref{pEigenvalue}, that
 $$ p_V(\dmat{\omega}^{-1} v \dmat{\omega})
 =  p_{g_2^{-1} U g_2}\bigl( g_2^{-1}( \dmat{\omega} u \dmat{\omega}^{-1} ) g_2
\bigr)
 = \omega^{2} \cdot p_{g_2^{-1} U g_2}( g_2^{-1} u g_2)
 = \omega^{2} \cdot p_V( v) .$$
 This establishes the second half of~\pref{ArchimedeanBddPf-pU=pVopp-omega}.

\begin{step} \label{ArchideanBddPf-contra}
 We obtain a contradiction.
 \end{step}
 Choose a unit $\omega \in \ints$, such that $\dmat{\omega} \in \Gamma$ and $\omega$ is not a root of unity. (This is possible because $\ints$  has infinitely many units.)
Step~\fullref{ArchimedeanBddPf-pU=pVopp}{omega}
implies that $\omega^2$ is real, so (by passing to a power) there is no harm in assuming $\omega >
1$.
 In the notation of Corollary~\fullref{OppScalarsSame}{scalars}, with $t = \dmat{\omega}$,
Step~\fullref{ArchimedeanBddPf-pU=pVopp}{omega}
asserts that
 $$ \mbox{$\omega_U = \omega^{-2} < 1$
 \quad and \quad
 $\omega_V = \omega^2 > 1$.}$$
 This contradicts Corollary~\fullref{OppScalarsSame}{same}, and thereby completes the proof of Theorem~\ref{ArchimedeanUBddOrbits}.
 \hfill\qedsymbol
 
 \begin{rem}
 For the lattices of $\Qrank$~one in $\SL(3,\real)$ that are discussed in the following section, we have $\omega_U = \omega_V$, instead of $\omega_U = 1/\omega_V$. For this reason, it is not so easy to obtain a contradiction for those groups.
 \end{rem}

\section{Lattices in $\SL(3,\real)$ or $\SL(3,\complex)$} \label{LattSL3NotRO}

In this section, we prove the following theorem.

\begin{thm} \label{GammaInSL3->UBdd}
Assume
\begin{itemize}
\item $G$ is either\/ $\SL(3,\real)$ or\/ $\SL(3,\complex)$,
\item $\Gamma$ is a noncocompact lattice in~$G$, 
\item we have a continuous action of\/~$\Gamma$ on\/~$\real$,
and
\item either $G = \SL(3,\real)$, or\/ $\Gamma$ does not contain any subgroup that is isomorphic to a noncocompact lattice in\/ $\SL(3,\real)$.
\end{itemize}
Then every orbit of every unipotent subgroup of\/~$\Gamma$ is a bounded subset of\/~$\real$.
\end{thm}

\begin{assump}
Throughout this section, $G$ and~$\Gamma$ are as described in Theorem~\ref{GammaInSL3->UBdd}.
\end{assump}

\begin{notation}
Write $G = \SL(3,\F_\infty)$, so $\F_\infty = \real$ or~$\complex$.
\end{notation}

The Margulis Arithmeticity Theorem provides a precise algebraic description of some finite-index subgroup of the lattice~$\Gamma$. Because there is no harm in replacing $\Gamma$ with this subgroup, we may assume the description applies to $\Gamma$ itself:

\begin{lem} \label{SL3existsKF}
We may assume there exist fields\/ $\F$ and\/~$\K$, such that
\begin{enumerate}
\item \label{SL3existsKF-K}
 $\K$ is a quadratic extension of\/~$\F$,
so\/ $\K = \F \bigl[ \sqrt{r} \bigr]$, for some $r \in \F$,
\item \label{SL3existsKF-either}
 either 
\begin{enumerate}
\item $\F_\infty = \real$, $\F = \rational$ and\/ $\K \subset \real$,
or
\item \label{SL3existsKF-either-C}
 $\F_\infty = \complex$, $\F$ is an imaginary quadratic extension of\/~$\rational$, and\/ $\K \cap \real = \rational$,
\end{enumerate}
and
\item $\Gamma$ is a finite-index subgroup of
  $$\SUJ(\ints) = \{\, g \in \SL(3,\ints) \mid g J \cjg{g}^\transpose
= J \,\} ,$$
 where
 \begin{itemize}
 \item $\ints$ is the ring of integers of\/~$\K$.
 \item $J = \begin{bmatrix}
 0 & 0 & 1 \\
 0 & 1 & 0 \\
 1 & 0 & 0
 \end{bmatrix}$,
 \item $\cjg{\phantom{g}}$ denotes the nontrivial Galois automorphism of the quadratic extension\/ $\K/\F$,
 and
 \item $\phantom{g}^\transpose$ denotes the transpose.
 \end{itemize}
\end{enumerate}
\end{lem}

\begin{proof}
Because $\Rrank G = 2 > 1$, the Margulis Arithmeticity Theorem \cite[Thm.~8.1.11, p.~298]{MargulisBook} tells us that $\Gamma$ is an arithmetic subgroup of~$G$. (Note that, since $\Gamma$ is not cocompact, we have no need to allow compact factors in the definition of an arithmetic subgroup \cite[Rem.~9.1.6(iii), p.~294]{MargulisBook}.) Thus, there is an algebraic number field~$\F$, with ring of integers~$\ints$, and an $\F$-form~$\G$ of~$G$, such that 
\begin{itemize}
\item either $\F_\infty = \real$ and $\F = \rational$, or $\F_\infty = \complex$ and $\F$ is an imaginary quadratic extension of~$\rational$,
and
\item (after passing to a finite-index subgroup) $\Gamma$ is isomorphic to a finite-index subgroup of the group of $\ints$-points of~$\G$.
\end{itemize}
Because $\Gamma$ is \emph{not} cocompact, we know $\Frank \G > 0$ (cf.\ \cite[Thm.~4.12, p.~210]{PlatonovRapinchukBook}). On the other hand, because  $\Gamma$ acts on $\real$, we must have $\Frank \G < 2$ \cite{Witte-ActOnCircle}. Therefore $\Frank \G = 1$.

Because $\Frank \G = 1$, the classification of $\F$-forms of $\SL(3,\F_\infty)$ \cite[Props.~2.17 and 2.18, pp.~87 and~88]{PlatonovRapinchukBook} asserts that $\Gamma$ must be exactly as described, except the requirement that $\K \cap \real = \rational$ when $\F_\infty = \complex$.

To complete the proof, suppose $\F_\infty = \complex$ and $\K \cap \real \neq \rational$. Then $\K \cap \real$ is a real quadratic extension of~$\rational$. Letting $\ints_\real = \ints \cap \real$ be the ring of integers of $\K \cap \real$, we see that $\SUJ(\ints \cap \real)$ is a noncocompact lattice in $\SL(3,\real)$. This contradicts the assumption that $\Gamma$ does not contain any noncocompact lattice in $\SL(3,\real)$.
\end{proof}

\begin{notation} \label{UinSL3Notn} \ 
 \begin{enumerate}

 \item For $\alpha,\zeta \in \K$ with $\zeta + \cjg\zeta = -\alpha \cjg\alpha$, let 
  $$ \mbox{$u(\alpha,\zeta) = \begin{bmatrix}
 1 & \alpha & \zeta \\
 0 & 1 & -\cjg\alpha \\
 0 & 0 & 1
 \end{bmatrix}$,
 \quad and \quad
 $v(\alpha,\zeta) = \begin{bmatrix}
 1 & 0 & 0 \\
 -\cjg\alpha & 1 & 0 \\
 \zeta  & \alpha & 1
 \end{bmatrix}$.}$$
 Note that $u(\alpha,\zeta)$ and $v(\alpha,\zeta)$ both belong to
$\SUJ(\K)$ (because $\zeta + \cjg\zeta = -\alpha \cjg\alpha$).

 \item Let
 $$ \mbox{$\U = \bigset{ u(\alpha, \zeta) }
 { \begin{matrix} \alpha, \zeta \in \K \\ \zeta + \cjg\zeta = -\alpha \cjg\alpha \end{matrix}  }$
 \quad and \quad
 $\V = \bigset{ v(\alpha,\zeta) }
 { \begin{matrix} \alpha, \zeta \in \K \\ \zeta + \cjg\zeta = -\alpha \cjg\alpha \end{matrix}} $},$$
 so $\U$ and~$\V$ are opposite maximal unipotent subgroups of
$\SUJ(\K)$. 

\item Let
 $$ \mbox{$U = \U \cap \Gamma$ \quad and \quad $V = \V \cap \Gamma$.}$$

\item For $\omega \in \ints$ with $\omega \cjg\omega = 1$, let
 $$ \dmat{\omega} = \begin{bmatrix}
 \omega & 0 & 0 \\
 0 & \cjg\omega^2 & 0 \\
 0 & 0 & \omega
 \end{bmatrix} .$$
 Note that  $\dmat{\omega} \in \SUJ(\ints)$.
 
 \item \label{UinSL3Notn-m}
  Because $\Gamma$ has finite index in $\SUJ(\ints)$, we may fix a positive
integer~$m$, such that 
 $$ \mbox{if $\alpha, \zeta \in m\ints$ (with $\zeta + \cjg\zeta = -\alpha \cjg\alpha$), 
 then $u(\alpha,\zeta) \in \Gamma$ and $v(\alpha,\zeta) \in \Gamma$.}$$

 \end{enumerate}

 \end{notation}
 
 The following calculation of Raghunathan is crucial when $\F = \rational$.

\begin{lem}[{(Raghunathan \cite[Lem.~1.7]{Raghunathan-CSP2})}] \label{RagId}
 Suppose $g$ and $\dmat{\omega}$ are elements of\/~$\Gamma$, with
 $$\mbox{$g = \begin{bmatrix}
 a & b & c \\
 *&*&* \\
 *&*&* 
 \end{bmatrix}$,
 \ 
 $\dmat{\omega} = \begin{bmatrix}
 \omega & 0 & 0 \\
 0 & \cjg\omega^2 & 0 \\
 0 & 0 & \omega
 \end{bmatrix}$,
 \ 
$\omega \cjg\omega = 1$,
 \  and 
 \  $\omega \equiv 1 \pmod{a \cjg a m}$.}$$
 If we let
 $\eta = (\omega^3 - 1)b/a$
 and
 $\xi = b \cjg\eta/a$,
 then $\bigl( g \, u(\eta,\xi) \bigr) (\dmat{\omega} g \dmat{\omega}^{-1})^{-1} \in
V$.
 \end{lem}

\begin{proof}
 Note that $\eta,\xi
\in m \ints$ (because $\omega \equiv 1 \pmod{a \cjg a m}$) and
 \begin{align*}
 -\eta \cjg \eta
 &= - \bigl( (\omega^3 - 1)b/a \bigr) \bigl( (\cjg \omega^3 - 1)\cjg b/ \cjg a \bigr) 
 \\&= -(\omega^3 - 1) (\cjg \omega^3 - 1) b \cjg b / ( a \cjg a)
 \\&= -(1 - \omega^3 - \cjg \omega^3 + 1) b \cjg b / ( a \cjg a)
  \\&= (\cjg \omega^3 - 1) b \cjg b / ( a \cjg a) + (\omega^3 - 1) b \cjg b / ( a \cjg a) 
  \\&= \cjg \eta b/a + \eta \cjg b / \cjg a
  \\&= \xi + \cjg \xi
  , \end{align*}
  so $u(\eta,\xi) \in \Gamma$. Therefore, 
 $$\bigl( g \, u(\eta,\xi) \bigr) (\dmat{\omega} g \dmat{\omega}^{-1})^{-1} \in
\Gamma .$$
 Easy calculations show that the two matrices $g \, u(\eta,\xi)$ and
$\dmat{\omega} g \dmat{\omega}^{-1}$ have the same first row, namely 
 $ \begin{bmatrix} a & \omega^3 b & c \end{bmatrix} $.
 Therefore, the first row of the product $\bigl(g \, u(\eta,\xi) \bigr)
(\dmat{\omega} g \dmat{\omega}^{-1})^{-1}$ is the same as the first row of
$(\dmat{\omega} g \dmat{\omega}^{-1}) (\dmat{\omega} g \dmat{\omega}^{-1})^{-1} = \Id$.
This means that the first row of the product is 
 $ \begin{bmatrix} 1 & 0 & 0 \end{bmatrix} $,
 so the product belongs to~$V$.
 \end{proof}
 
 Recall that $r$ is an element of~$\F$, such that $\K = \F \bigl[ \sqrt{r} \bigr]$ \fullsee{SL3existsKF}{K}.

\begin{cor} \label{RagIdPrecise}
 Given $\beta \in 2m \ints$, $\ell \in m \integer$, and $\omega \in \ints$,
such that
 $$ \mbox{$\omega \cjg\omega = 1$
 \quad and \quad
 $\omega \equiv 1 \bigpmod{(1 - \ell^2 y^2 r)m}$,
 where $y = -\beta \cjg\beta / 2$} ,$$
 let 
 $$
 \eta = \frac{(\omega^3 - 1) \ell \beta \sqrt{r}}{1 + \ell y
\sqrt{r}} 
 \quad and \quad
 \lambda = \frac{(1 - \omega^3) \beta}{1 - y \ell \sqrt{r}}
 . $$
 Then there exist 
 $$\mbox{
 $v(\lambda) = v(\lambda,*)$,
 $v(\beta) = v(\beta,*)$,
 $v(\omega^3 \beta) = v(\omega^3 \beta,*)$,
 }$$
 $$
 \mbox{
 $u(\eta) = u(\eta,*)$,
 and
 $z = u(0, \ell \sqrt{r})$
 }$$
 in~$\Gamma$, such that
 $$ v(\beta) \,u(\eta) = \bigl( z^{-1} \, v(\lambda) \, z \bigr) \,
v(\omega^3 \beta) .$$
 \end{cor}

\begin{proof}
 Let
 \begin{itemize}
 \item $v(\beta) = v(\beta,y)$ and
 $v(\omega^3 \beta) = \dmat{\omega} \, v(\beta) \, \dmat{\omega}^{-1} =
v(\omega^3 \beta, y)$,
 and

 \item $u(\eta) = u(\eta,\xi)$, where 
 $$\xi = \frac{\ell \beta \cjg\eta \sqrt{r}}{1 + \ell y \sqrt{r}} 
 = \frac{ \eta \cjg\eta }{ \omega^3 - 1} .$$
 \end{itemize}
 Note that $v(\beta)$, $v(\omega^3 \beta)$, $u(\eta)$, and~$z$ are
elements of $\Gamma$ \fullcf{UinSL3Notn}{m}.

By letting 
 $$g = z \, v(\beta)
 = \begin{bmatrix}
  1&0&\ell \sqrt{r}\\
  0&1&0 \\
 0&0&1 
 \end{bmatrix}
 \begin{bmatrix}
  1&0&0 \\
  -\cjg\beta&1&0 \\
 y&\beta&1 
 \end{bmatrix}
 = \begin{bmatrix}
  1 + \ell y \sqrt{r} & \ell \beta \sqrt{r} & \ell \sqrt{r}\\
  *&*&* \\
 *&*&* 
 \end{bmatrix},$$
 we see, from Lemma~\ref{RagId}, that there exists $v \in V$,
 such that 
 $z \, v(\beta) \, u(\eta) =  v \,\bigl( \dmat{\omega} \, z \, v(\beta) \,
\dmat{\omega}^{-1} \bigr)$.
 Because 
 \begin{itemize}
 \item $z$ commutes with $\dmat{\omega}$, 
 and
 \item $\dmat{\omega} \, v(\beta) \, \dmat{\omega}^{-1} = v(\omega^3 \beta)$,
 \end{itemize}
 we conclude that
 $$ v(\beta) \, u(\eta) =  (z^{-1} v z) \, v(\omega^3 \beta) .$$
 Writing 
 $$v = v(\mu, *) ,$$
 we see that:
 \begin{itemize}
 \item the $(2,1)$ entry of $v(\beta) \, u(\eta)$ is
 $$ (-\cjg\beta)(1) + (1)(0) +(0)(0) = -\cjg\beta ,$$
 \item the second row of $z^{-1} v z$ is
 $\begin{bmatrix} -\cjg\mu & 1 & -\ell \cjg\mu \sqrt{r}
\end{bmatrix}$,
 \item the $(2,1)$ entry of $(z^{-1} v z) \, v(\omega^3 \beta)$ is
 $$ (-\cjg\mu)(1) + (1) (-\cjg\omega^3 \cjg\beta) + (-\ell
\cjg\mu \sqrt{r}) (y)
 = -\cjg\mu ( 1 + y \ell \sqrt{r} ) -\cjg\omega^3 \cjg\beta .$$
 \end{itemize}
 The $(2,1)$ entries that we calculated must be equal, so we conclude that 
 $$ \mu =  \frac{(1 - \omega^3) \beta}{1 - y \ell \sqrt{r}} =
\lambda.$$
 Thus, we may let $v(\lambda) = v$.
 \end{proof}

\begin{proof}[{\bf Proof of Theorem~\ref{GammaInSL3->UBdd}}]
 Assume some orbit of some unipotent subgroup $U_0$ is not bounded. (This will lead to a contradiction.) There is no harm in assuming, for definiteness, that:
 \begin{enumerate}
 \item $U_0 = U$ is a maximal unipotent subgroup,
 and
 \item the $U$-orbit of~$0$ is not bounded above.
 \end{enumerate}
 
 The proof now proceeds in a sequence of steps.

\setcounter{step}{0}

\begin{step}
 We may assume that the $V$-orbit of\/~$0$ is not bounded
above.
 \end{step}
 See the argument at the start of Step~\ref{ArchimedeanBddPf-pU=pVopp} of the proof of Theorem~\ref{ArchimedeanUBddOrbits} in \S\ref{ArchimedeanUBddOrbitsPf}.

\begin{step} \label{GammaInSL3->UBddPf-pU}
 We may assume
 \begin{enumerate} \renewcommand{\theenumi}{\alph{enumi}}
 \item \label{GammaInSL3->UBddPf-pU-Finfty=R}
  $\F_\infty = \real$,
 and
 \item \label{GammaInSL3->UBddPf-pU-pU}
 $p_U \bigl( u(\alpha,\zeta) \bigr) = \alpha$, for all
$u(\alpha,\zeta) \in U$.
\end{enumerate}
 \end{step}
 It is well known (cf.~\cite[Prop.~2.15(3), p.~86]{PlatonovRapinchukBook}) that $\SL(3,\F_\infty)$ can be viewed as the $\real$-points of an almost simple algebraic $\rational$-group, such that 
 \begin{itemize}
 \item $\Gamma$ is an arithmetic subgroup of $\SL(3,\F_\infty)$,
 and
 \item $\GU = \begin{bmatrix} 1&*&*\\ 0&1&*\\ 0&0&1 \end{bmatrix}$
 and
 $\GT = \begin{bmatrix} *&0&0\\ 0&*&0\\ 0&0&* \end{bmatrix}$
 are $\rational$-subgroups of $\SL(3,\F_\infty)$.
 \end{itemize}
 
 \pref{GammaInSL3->UBddPf-pU-Finfty=R}
 Suppose $\F_\infty = \complex$. Then the group $\ints^\times \cap \real$ of real units of~$\ints$ is finite \fullfullcf{SL3existsKF}{either}{C}, so no finite-index subgroup of~$\ints^\times$ is contained in~$\real$. Therefore, no one-parameter subgroup of~$\GU$ is normalized by $\GT \cap \Gamma$. This contradicts the existence of the subgroup~$\GU_1$ provided by Corollary~\ref{p=root}. 
 
 \pref{GammaInSL3->UBddPf-pU-pU}
 We may now assume $\F_\infty = \real$. The only one-parameter subgroups of~$\GU$ that are normalized by $\GT \cap
\Gamma$ are the root subgroups
 $$  \mbox{$\GU_\alpha = \begin{bmatrix} 1&*&0\\ 0&1&0\\ 0&0&1 \end{bmatrix}$,
 \quad
 $\GU_{\cjg\alpha} = \begin{bmatrix} 1&0&0\\ 0&1&*\\ 0&0&1 \end{bmatrix}$,
 \quad and \quad
 $\GU_\zeta = \begin{bmatrix} 1&0&*\\ 0&1&0\\ 0&0&1 \end{bmatrix}$.}$$
 Now $\GU_\zeta = [\GU,\GU]$ is in the kernel of~$p_U$, so we conclude that the
one-parameter subgroup of Corollary~\ref{p=root} is either $\GU_\alpha$ or~$\GU_{\cjg\alpha}$.
Thus, up to a (nonzero) scalar multiple, $p_U \bigl( u(\alpha,*) \bigr)$ is
either $\alpha$ or~$\cjg\alpha$. This implies that, up to a \emph{positive} scalar multiple,
$p_U \bigl( u(\alpha,*) \bigr)$ is either $\alpha$, $-\alpha$, $\cjg\alpha$,
or~$-\cjg\alpha$. Because 
 \begin{itemize}
 \item replacing $\Gamma$ with its Galois conjugate
$\cjg\Gamma$ transforms $u(\alpha,*)$ to
$u(\cjg\alpha,*)$, 
 and
 \item conjugation by the matrix $\dmat{\omega}$ with $\omega = -1$ transforms $u(\alpha,*)$ to
$u(-\alpha,*)$, 
 \end{itemize}
 we may assume that $p_U \bigl( u(\alpha,*) \bigr) = \alpha$.

\begin{step} \label{SUPf-pV=beta}
 We may assume $p_V \bigl( v(\beta,\xi) \bigr) = \pm \beta$, for all
$v(\beta,\xi) \in V$.
 \end{step}
 Suppose $p_V \bigl( v(\beta,\xi) \bigr) = \pm \cjg\beta$. (This will lead to a contradiction.) Choose a unit $\omega \in \ints$, such that $\dmat{\omega} \in \Gamma$ and $\omega$ is not a root of unity. (This is possible because $\ints$  has infinitely many units.)
Step~\fullref{GammaInSL3->UBddPf-pU}{Finfty=R}
implies that $\omega$ is real, so (by replacing $\omega$ with~$\omega^{\pm 2}$) there is no harm in assuming $\omega > 1$ and $\cjg\omega > 0$.
 In the notation of Corollary~\fullref{OppScalarsSame}{scalars}, with $t = \dmat{\omega}$,
 Step~\fullref{GammaInSL3->UBddPf-pU}{pU} implies that 
 $$ \omega_U = \omega^{-3} < 1 .$$
Our assumption that $p_V \bigl( v(\beta,\xi) \bigr) = \pm \cjg\beta$ implies
 $$ \omega_V = \cjg\omega^{-3} = \omega^3 > 1 .$$
 This contradicts \fullref{OppScalarsSame}{same}.

(Alternatively, this step may be justified by the argument in
Step~\fullref{ArchimedeanBddPf-pU=pVopp}{omega} on page~\pageref{ArchimedeanBddPf-pU=pVopp-omega}.
In fact, this method yields the more precise result that 
 $p_V \bigl( v(\beta,*) \bigr) = - \beta$.) 
 
\begin{step}
 Fix some nonzero $\ell \in m \integer$, and let $z = u(0,\ell \sqrt{r})$; we may assume
that $z$ fixes\/~$0$.
 \end{step}
 Because $z \in [\U,\U]$, we know that $p_U(z) = 0$. Therefore,
Theorem~\fullref{AultRhemtulla}{p=0} implies that $z$ has a fixed point~$x$ in
the interval $[0, 0 \cdot u]$ (for any $u \in U$ with $p(u) > 0$). There is
no harm in assuming $x = 0$.

\begin{step} \label{SUPfMainIneq}
 There exist nonzero $\alpha,\beta \in 2m \ints$ and $\omega \in \ints$, such
that
 \begin{enumerate} \renewcommand{\theenumi}{\alph{enumi}}
 \item \label{SUPfMainIneq-ineq}
 $0 \cdot v(\beta,*) < 0 \cdot u(\alpha,*) < 0 \cdot v(\omega^3 \beta,*)$,
 \item \label{SUPfMainIneq-omega}
 $\omega \cjg\omega = 1$, $\omega > 1$, $\dmat{\omega} \in \Gamma$, and
$\omega \equiv 1 \bigpmod{(1 - \ell^2 y^2 r)m}$,
 where $y = - \beta \cjg\beta / 2$,
 \item \label{SUPfMainIneq-pV}
 $p_V \bigl( v(\beta,*) \bigr) > 0$,
 and
 \item \label{SUPfMainIneq-ConjBeta}
 $\cjg\beta > 0$.
 \end{enumerate}
 \end{step}
 Because the sign of~$\cjg\beta$ is independent of the sign of~$\beta$ (and using Step~\ref{SUPf-pV=beta}),
there exists $\beta \in 2m \ints$, satisfying \pref{SUPfMainIneq-pV}
and~\pref{SUPfMainIneq-ConjBeta}.

 Because the $U$-orbit of~$0$ is not
bounded above (but each orbit of $[U,U]$ or $[V,V]$ is bounded
\fullcf{AultRhemtulla}{p=0}), there exists $\alpha \in 2m \ints$, such that
$0 \cdot v(\beta,*) < 0 \cdot u(\alpha,*)$. 
 Similarly, because the $V$-orbit of~$0$ is not bounded above, 
 there exists $\omega \in \ints$, such that $\omega \cjg\omega = 1$, $\omega > 1$, and $0 \cdot u(\alpha,*) < 0 \cdot v(\omega^3 \beta,*)$. This
establishes~\pref{SUPfMainIneq-ineq}.

By replacing $\omega$ with an appropriate power~$\omega^n$ (with $n > 0$), we
obtain the conditions of~\pref{SUPfMainIneq-omega}.

\begin{step}
 We obtain a contradiction.
 \end{step}
 Letting
 $$ 
 \eta = \frac{(\omega^3 - 1) \ell \beta \sqrt{r}}{1 + \ell y \sqrt{r}} 
 \mbox{\quad and \quad}
 \lambda = \frac{(1 - \omega^3) \beta}{1 - y \ell \sqrt{r}}
 , $$
 we know, from Corollary~\ref{RagIdPrecise}, that there exist $v(\lambda)$, $v(\beta)$,  $v(\omega^3 \beta)$,
$u(\eta)$, and $z = u(0,\ell \sqrt{r})$, such that
 \begin{equation} \label{LucyEq}
 v(\beta) \,u(\eta) = \bigl( z^{-1} \,
v(\lambda) \, z \bigr) \, v(\omega^3 \beta) .
 \end{equation}
 From Step~\fullref{SUPfMainIneq}{ineq}, we know $0 \cdot v(\beta)
< 0 \cdot u(\alpha)$. 
Therefore
 \begin{equation} \label{SUPf-vu<uu}
 0 \cdot v(\beta) \, u(\eta) \< 0 \cdot u(\alpha) \, u(\eta) .
 \end{equation}
 We may assume $\ell$ is large enough that $|\ell \sqrt{r}| > 2$, so $1 + \ell y \sqrt{r}$ has the same sign as $\ell y \sqrt{r}$. Then, because 
 $$ \mbox{$\omega^3 - 1 > 0$, 
 \quad 
 $y = - \beta \cjg\beta/2$, 
 \quad and \quad
 $\cjg\beta > 0$,} $$
 we see, from the definition of~$\eta$, that $\eta < 0$. Therefore
 $$ p_U \bigl( u(\alpha) \, u(\eta) \bigr)
 = p_U \bigl( u(\alpha + \eta, *) \bigr)
 = \alpha + \eta
 < \alpha
 = p_U \bigl( u(\alpha) \bigr) ,$$
 so
 \begin{equation} \label{SUPf-uu<u}
 0 \cdot u(\alpha) \, u(\eta) \< 0 \cdot u(\alpha) .
 \end{equation}
  From Step~\fullref{SUPfMainIneq}{ineq}, we have
 \begin{equation} \label{SUPf-u<v}
 0 \cdot u(\alpha) \< 0 \cdot v(\omega^3 \beta) .
 \end{equation}
  Replacing $\ell$ with~$-\ell$ would not require any change in $\alpha$, $\beta$, or~$\omega$ (because the conditions in Step~\fullref{SUPfMainIneq}{omega} depend only on~$\ell^2$, not on~$\ell$). Thus, we may assume $p_V \bigl( v(\lambda) \bigr) > 0$ \cfStep{SUPf-pV=beta}. Then $0 \< 0 \cdot v(\lambda)$.
Because $z$ fixes~$0$, this implies that $0 \< 0 \cdot z^{-1} \, v(\lambda) \, z$,
so
 \begin{equation} \label{SUPf-v<vv}
 0 \cdot v(\omega^3 \beta) \< 0 \cdot \bigl( z^{-1}
\, v(\lambda) \, z \bigr) \, v(\omega^3 \beta)  .
 \end{equation}
 By combining \pref{SUPf-vu<uu}, \pref{SUPf-uu<u}, \pref{SUPf-u<v}, and
\pref{SUPf-v<vv}, we conclude that
 $$ 0 \cdot v(\beta) \, u(\eta) \< 0 \cdot \bigl( z^{-1}
\, v(\lambda) \, z \bigr) \, v(\omega^3 \beta) .$$
 This contradicts \pref{LucyEq}.
 \end{proof}

\section{Proof of Theorem~\ref{AllNotRO}} \label{AllLattSect}

 Throughout this section, the conditions in the statement of Theorem~\ref{AllNotRO} are assumed to be satisfied:
\begin{itemize}
\item Conjecture~\ref{SU2conj} is true,
\item $G$ is a connected, semisimple Lie group with finite center,
\item $\Rrank G \ge 2$,
and
\item $\Gamma$ is a noncocompact, irreducible lattice in~$G$.
\end{itemize}
We wish to show that $\Gamma$ has no nontrivial, orientation-preserving action on~$\real$.

Because $\real$ has no orientation-preserving homeomorphisms of finite order, there is no harm in modding out a finite group. Thus, we may assume that $G$ has trivial center. Hence, $G$ is linear (and Corollary~\ref{ActMustBeFaithful} implies that every nontrivial, orientation-preserving action of~$\Gamma$ on~$\real$ is faithful). We now recall the following theorem:

 \begin{thm}[{(Chernousov-Lifschitz-Morris \cite{CLMmin})}] \label{CLMminVague}
 Assume
 \begin{itemize}
 \item $G$ is a connected, semisimple, linear Lie group,
 \item $\Rrank G \ge 2$,
 and
 \item $\Gamma$ is a noncocompact, irreducible lattice in~$G$.
 \end{itemize}
Then\/ $\Gamma$ contains a subgroup that is isomorphic to either
\begin{enumerate}
\item \label{CLMminVague-SL3}
 a noncocompact lattice in\/ $\SL(3,\real)$ or\/ $\SL(3,\complex)$,
or
\item \label{CLMminVague-SL2}
 a finite-index subgroup of\/ $\SL(2,\ints)$, where
\begin{itemize}
\item $\ints$ is the ring of integers of a number field\/~$\F$,
and
\item $\F$ is neither\/ $\rational$ nor an imaginary quadratic extension of\/~$\rational$.
\end{itemize}
\end{enumerate}
 \end{thm}

By passing to a subgroup, we may assume $\Gamma$ is as described in either~\fullref{CLMminVague}{SL3} or~\fullref{CLMminVague}{SL2}. Theorem~\ref{SL(2,O)NoActR} tells us that the groups in~\fullref{CLMminVague}{SL2} have no nontrivial, orientation-preserving actions on~$\real$, so we may assume $\Gamma$ is a noncocompact lattice in either $\SL(3,\real)$ or $\SL(3,\complex)$. 

Furthermore, by passing to a subgroup again, we may assume that either 
\begin{itemize}
\item $\Gamma$ is a noncocompact lattice in $\SL(3,\real)$,
 or
 \item $\Gamma$ does not contain any subgroup that is isomorphic to a noncocompact lattice in $\SL(3,\real)$. 
\end{itemize}
Now Theorem~\ref{GammaInSL3->UBdd} tells us, for every orientation-preserving action of~$\Gamma$ on~$\real$, that every orbit of every unipotent subgroup of\/~$\Gamma$ is a bounded subset of\/~$\real$. Also, because Conjecture~\ref{SU2conj} is assumed to be true, we know that $\Gamma$ is virtually boundedly generated by unipotents. So Lemma~\ref{BddGen+BddOrbits} implies that $\Gamma$ has no nontrivial, orientation-preserving action on~$\real$. This completes the proof of Theorem~\ref{AllNotRO}.
\hfill \qedsymbol

\end{document}